\documentclass[a4paper,12pt,final]{amsart}
\usepackage{times,a4wide,mathrsfs,amssymb,amsmath,amsthm,enumerate,xypic,tikzsymbols,dsfont}

\newcommand{\C}{\mathbb{C}}
\newcommand{\ZZ}{\mathbb{Z}}

\newcommand{\QQ}{\mathbb{Q}}
\newcommand{\NN}{\mathbb{N}}
\newcommand{\PP}{\mathbb{P}}

\newcommand{\MM}{\mathcal M}

\DeclareMathOperator{\aut}{Aut}

\DeclareMathOperator{\sym}{Sym}

\DeclareMathOperator{\jac}{Jac}

\newtheorem{theorem}{Theorem}[section]

\newtheorem{corollary}[theorem]{Corollary}
\newtheorem{proposition}[theorem]{Proposition}

\newtheorem{convention}{Conventions}

\newtheorem{nonumbering}{Theorem}

\theoremstyle{definition}
\newtheorem{remark}[theorem]{Remark}
\newtheorem{definition}[theorem]{Definition}

\newtheorem{notation}[theorem]{Notation}

\newtheorem{nonumberingt}{Acknowledgments}

\newtheorem{nonumberingda}{Data availability statement}

\newtheorem{nonumberingcon}{Conflict of interest statement}

\begin{document}

\author[Robert Laterveer]
{Robert Laterveer}

\address{Institut de Recherche Math\'ematique Avanc\'ee,
CNRS -- Universit\'e 
de Strasbourg,\
7 Rue Ren\'e Des\-car\-tes, 67084 Strasbourg CEDEX,
FRANCE.}
\email{robert.laterveer@math.unistra.fr}

\title{On the tautological ring of Humbert curves}

\begin{abstract} We exhibit a 2-dimensional family of non-hyperelliptic curves of genus 5, called Humbert curves, for which the tautological ring injects into cohomology. In particular, Humbert curves have a multiplicative Chow--K\"unneth decomposition (in the sense of Shen--Vial), and their Ceresa cycle is torsion.
 \end{abstract}

\thanks{Supported by ANR grant ANR-20-CE40-0023.}

\maketitle

\section{Introduction}

Given a smooth projective variety $Y$ over $\C$, let $A^i(Y):=CH^i(Y)_{\QQ}$ denote the Chow groups of $Y$ (i.e. the groups of codimension $i$ algebraic cycles on $Y$ with $\QQ$-coefficients, modulo rational equivalence \cite{F}, \cite{Vo}). The intersection product defines a ring structure on $A^\ast(Y)=\bigoplus_i A^i(Y)$, the {\em Chow ring\/} of $Y$.

Motivated by the particular behaviour of Chow rings of K3 surfaces \cite{BV} and of abelian varieties \cite{Beau}, Beauville \cite{Beau3} has famously conjectured that for certain special varieties, the Chow ring should admit a multiplicative splitting. To make concrete sense of Beauville's elusive ``splitting property'' conjecture, Shen--Vial \cite{SV} introduced the concept of {\em multiplicative Chow--K\"unneth decomposition\/} (we will abbreviate to ``MCK decomposition''). In short, this is a graded decomposition of the Chow motive of a smooth projective variety, such that the intersection product respects the grading (cf. subsection \ref{ss:mck} for details).

It is something of a challenge to understand the class of varieties admitting an MCK decomposition: abelian varieties, K3 surfaces and cubic hypersurfaces are in this class. As for curves, it is known that hyperelliptic curves are in this class, whereas not all curves are in this class.

The main result of the present paper is about {\em Humbert curves\/}; these are certain non-hyperelliptic curves of genus 5 that were studied as far back as 1894 (cf. \cite{Hum} and subsection \ref{ssh} below).

\begin{nonumbering}[=Theorems \ref{main} and \ref{main1}] Let $C$ be a Humbert curve, and $m\in\NN$. Let
  \[ R^\ast(C^m):= \bigl\langle  (p_{ij})^\ast(\Delta_C), (p_k)^\ast(K_C)\bigr\rangle \ \ \subset\ A^\ast(C^m) \]
  be the $\QQ$-subalgebra generated by (pullbacks of) the diagonal $\Delta_C\subset C\times C$ and (pullbacks of) the canonical divisor $K_C$. The cycle class map induces injections
   \[ R^\ast(C^m)\ \hookrightarrow\ H^\ast(C^m,\QQ) \]
   for all $m\in\NN$.
   
 In particular, $C$ has an MCK decomposition, and the Faber--Pandharipande cycle
   \[    \Delta_C\cdot (p_j)^\ast(K_C) - {1\over 8} K_C\times K_C   \ \ \in\  A^2(C\times C) \ \ \ \ (j=1,2)\]
   and the Ceresa cycle
    \[   [C] - (-1_{\jac(C)})_\ast [C]     \ \   \in\ A_1(\jac(C)) \]
    (for an appropriate choice of embedding of $C$ in its Jacobian $\jac(C)$) 
    are zero.  
\end{nonumbering}

The statement about the {\em tautological ring\/} $R^\ast(C^m)$ is inspired by Tavakol's work \cite{Ta2}, \cite{Ta}, where the same injectivity statement is proven for hyperelliptic curves.

Theorem \ref{main} provides the first examples of non-hyperelliptic curves with an MCK decomposition. 
These are also the first examples of non-hyperelliptic curves with torsion Ceresa cycle (non-hyperelliptic examples where the class of the Ceresa cycle in the intermediate Jacobian is torsion are given in \cite{BLLS} and \cite{Bcurve}, cf. also \cite{Bcurve2} for a result modulo algebraic equivalence\footnote{(Added in revision) Cf. also the recent work \cite{QZ}.}). The Ceresa cycle is famously known to be torsion for hyperelliptic curves, but non-torsion in general (cf. \cite{Cere} and subsection \ref{ssc} below). As for the Faber--Pandharipande cycle, this is known to be torsion for plane curves and for hyperelliptic curves, but non-torsion in general (cf. \cite{GG} and Remark \ref{!} below).

 \vskip0.5cm

\begin{convention} In this article, the word {\sl variety\/} will refer to a reduced irreducible scheme of finite type over $\C$. A {\sl subvariety\/} is a (possibly reducible) reduced subscheme which is equidimensional. 

{\bf All Chow groups are with rational coefficients}: we will denote by $A_j(Y)$ the Chow group of $j$-dimensional cycles on $Y$ with $\QQ$-coefficients; for $Y$ smooth of dimension $n$ the notations $A_j(Y)$ and $A^{n-j}(Y)$ are used interchangeably. 
%The notations $A^j_{hom}(Y)$ and $A^j_{AJ}(Y)$ will be used to indicate the subgroup of homologically trivial (resp. Abel--Jacobi trivial) cycles.
For a morphism $f\colon X\to Y$, we will write $\Gamma_f\in A_\ast(X\times Y)$ for the graph of $f$.

The contravariant category of Chow motives (i.e., pure motives with respect to rational equivalence as in \cite{Sc}, \cite{MNP}) will be denoted 
$\MM_{\rm rat}$.
\end{convention}

\section{Preliminaries}

\subsection{Humbert curves} 
\label{ssh}

\begin{definition} A smooth projective curve $C$ is called a {\em Humbert--Edge curve of type $n$\/} if $C\subset\PP^n$ can be defined as a complete intersection of $n-1$ quadrics
      \[   \begin{cases}  x_0^2+x_1^2+\cdots \  \cdots \ \cdots+ x_{n}^2&=0\ ,\\
                          \lambda_0 x_0^2+\lambda_1 x_1^2+\cdots + \lambda_{n}x_{n}^2&=0\ ,\\
                            \lambda_0^2 x_0^2+\lambda_1^2 x_1^2+\cdots + \lambda_{n}^2 x_{n}^2&=0\ ,\\
                          \ \ \ \vdots\\
                         \ \ \  \vdots\\
                           \lambda_0^{n-2} x_0^2+\lambda_1^{n-2} x_1^2+\cdots + \lambda_{n}^{n-2} x_{n}^2&=0\ ,\\     
                           \end{cases}\]
where $\lambda_i\in\C$, $i=0,\ldots,n$, are pairwise distinct complex numbers.

A {\em Humbert curve\/} is a Humbert--Edge curve of type $4$.
\end{definition}

\begin{theorem} A Humbert--Edge curve of type $n\ge 4$ is not hyperelliptic, and not trigonal.
\end{theorem}

\begin{proof} This is \cite[Theorems 2.1 and 2.2]{CGH}.
\end{proof}

 \begin{remark} A Humbert--Edge curve $C$ of type $n$ has genus 
   \[ g(C)=2^{n-2}(n-3)+1\ .\] 
   In particular, Humbert curves are of genus 5, and Humbert--Edge curves of type 2 and type 3 are quadrics resp. elliptic curves.
 
% It is known that Humbert--Edge curves of type $n\ge 4$ are not hyperelliptic, and not trigonal \cite[Theorems 2.1 and 2.2]{CGH}.
 
 An alternative coordinate-free definition is as follows: a smooth curve $C$ of genus $g(C)=2^{n-2}(n-3)+1$ is a Humbert--Edge curve of type $n$ if and only if $C$ admits an action of a group $G\cong(\ZZ/2\ZZ)^n$ with generators $\sigma_0,\ldots,\sigma_n$ such that $\sigma_0\cdots\sigma_n=1$ and the fixed loci of $\sigma_i$ are pairwise disjoint of cardinality $2^{n-1}$
 \cite{CGH}, \cite{FZ}.
 
 The Jacobian of Humbert--Edge curves is studied in \cite{FZ} and \cite{AAR}.
 
Humbert curves are named after Humbert \cite{Hum}, who first described them as certain curves in $\PP^3$: a Humbert curve is the locus of points of tangency of lines through a fixed point in $\PP^3$ and tangent to twisted cubics through 5 general points. Subsequently, this was generalized by Edge \cite{Edge}, \cite{Edge2}. Yet another characterization is as follows: a genus 5 curve is a Humbert curve if and only if it has 5 bi-elliptic involutions (i.e., involutions yielding a quotient curve of genus 1) \cite[Sections 2.1 and 2.2]{KMV}.
A remarkable relation between Humbert curves and Weddle surfaces via Prym varieties is exhibited in \cite{Var}.
 
 The moduli space of Humbert curves is 2-dimensional and is explicitly described in \cite[Section 2]{HR}. For $n\ge 5$ odd, the moduli space of Humbert--Edge curves of type $n$ is isomorphic to the moduli space of hyperelliptic curves of genus ${n-1\over 2}$ \cite[Section 3.3]{CF}.
  \end{remark}
 
 \begin{remark} As explained and exploited in \cite[Section 2]{AAR}, Humbert--Edge curves have a nice inductive structure. Given a Humbert--Edge curve $C\subset\PP^n$ of type $n$, let
 $\sigma_i$ ($i=0,\ldots,n$) be the automorphism of $C$ obtained by sending $x_i$ to $-x_i$. Then the quotient $C/\langle\sigma_i\rangle$ is a Humbert--Edge curve of type $n-1$.
 \end{remark}

 \subsection{The motive of Humbert curves} 
 
 \begin{notation} Let $C\subset\PP^4$ be a Humbert curve, defined by equations
     \[   \begin{cases}  x_0^2+x_1^2+\cdots \  \cdots \ \cdots+ x_{4}^2&=0\ ,\\
                          \lambda_0 x_0^2+\lambda_1 x_1^2+\cdots + \lambda_{4}x_{4}^2&=0\ ,\\
                            \lambda_0^2 x_0^2+\lambda_1^2 x_1^2+\cdots + \lambda_{4}^2 x_{4}^2&=0\ ,\\
                           \end{cases}\]
where $\lambda_0, \ldots,\lambda_4$, are pairwise distinct complex numbers.

For $i=0,\ldots,4$, let $\sigma_i\in\aut(C)$ be the involution of $C$ obtained by sending $x_i$ to $-x_i$. Let $G\cong(\ZZ/2\ZZ)^4\subset\aut(C)$ be the group with generators $\sigma_0,\ldots,\sigma_4$ and relation $\sigma_0\cdots\sigma_4=1$.
    \end{notation}

 \begin{proposition}\label{splitting} Let $C\subset\PP^4$ be a Humbert curve. There is a splitting of the Chow motive
   \[ h(C)= h^0(C) \oplus h^2(C) \oplus \bigoplus_{i=0}^4 h^1_{(i)}(C )\ \ \ \hbox{in}\ \MM_{\rm rat}\ ,\]
   with the property that $H^\ast(h^0(C),\QQ)=H^0(C,\QQ)$, $H^\ast(h^2(C),\QQ)=H^2(C,\QQ)$, and
     \[ \begin{split} H^\ast(h^1_{(i)}(C),\QQ) = \bigl\{  a\in H^1(C,\QQ)\,\vert \,(\sigma_i)_\ast(a)=a\ \hbox{and}\ (\sigma_j)_\ast(a)=-a\forall j\not=i\bigr\}&
     \cr   (i=0,\ldots,&4) \ .\cr
     \end{split} \]
     
     Moreover, writing $\rho_i\colon C\to C_i:=C/\langle\sigma_i\rangle$ for the quotient morphism, there are isomorphisms of Chow motives
     \[ {}^t \Gamma_{\rho_i}\colon\  h(C_i)\cong h^0(C)\oplus h^1_{(i)}(C)\oplus h^2(C)\ \ \ \hbox{in}\ \MM_{\rm rat}\ \ \ (i=0,\ldots,4)\ ,\]
     and the $C_i$ are elliptic curves.
      \end{proposition}  
    
    \begin{proof} Let $z\in A^1(C)$ be the unique 0-cycle of degree 1 that is $\sigma_i$-invariant for all $i$. Concretely, this 0-cycle can be constructed as follows: let $\rho\colon C\to C/G\cong \PP^1$ be the quotient morphism. Then $\rho$ has degree $2^4$ and one defines
    \[    z:= {1\over 2^4}\, \rho^\ast(p)\ \ \in\ A^1(C)\ ,\] 
    where $p\in C/G$ is an arbitrary point. This choice of 0-cycle defines a decomposition
     \[ h(C)= h^0(C)+ h^1(C) + h^2(C)\ \ \ \hbox{in}\ \MM_{\rm rat}\ ,\]
     where $h^0(C)$ and $h^2(C)$ are defined by the projectors $\pi^0_C:=z\times C$ resp. $\pi^2_C:=C\times z$.
     
%     Next, let us define projectors
%     \[  \pi^1_{C,(i)}:= {1\over 2^4}\ (\Delta_C+\Gamma_{\sigma_i})\circ \pi^1_C\ \ \in\ A^1(C\times C)\ \ \ (i=0,\ldots,4)\ \]
%     (these are projectors because of the commutativity $\Gamma_{\sigma_i}\circ \pi^k_C=\pi^k_C\circ\Gamma_{\sigma_i}$ for $k=0,2$ and hence also for $k=1$).
%    We note that $\pi^1_{C,(i)}$ is a projector on the $\sigma_i$-invariant part of $H^1(C,\QQ)$. That is, writing 

Let us write $\rho_i\colon C\to C_i:= C/\langle\sigma_i\rangle$ for the quotient morphism, so that $(\rho_i)^\ast H^1(C_i,\QQ)$ is the $\sigma_i$-invariant part.
  %  \[  (\pi^1_{C,(i)})_\ast H^\ast(C,\QQ)= (\rho_i)^\ast H^1(C_i,\QQ) \ .\]
    The curve $C_i$ is (a Humbert--Edge curve of type 3 and so) an elliptic curve. Moreover, for $i\not=j$ the intersection $ (\rho_i^\ast)^\ast H^1(C_i,\QQ) \cap  (\rho_j^\ast)^\ast H^1(C_j,\QQ)   $ is empty, since the quotient $C/\langle \sigma_i,\sigma_j\rangle$ is (a Humbert--Edge curve of type 2 and so) a rational curve. Since $\dim H^1(C,\QQ)=10$ and $\dim H^1(C_i,\QQ)=2$, it follows that  
  there is a decomposition of cohomology
    \[ H^1(C,\QQ)= \bigoplus_{i=0}^4  (\rho_i)^\ast H^1(C_i,\QQ) \ .\]
    To upgrade this decomposition to a motivic decomposition, one writes
    \[ \Delta_C={1\over 2} (\Delta_C+\Gamma_{\sigma_i}+\Delta_C -\Gamma_{\sigma_i})\ \ \ (i=0,\ldots,4)\ ,\]
    which gives
      \begin{equation}\label{this} \begin{split} \pi^1_C = {1\over 2^5}\, \sum_{r_0=0}^1 \sum_{r_1=0}^1 \cdots \sum_{r_4=0}^1 \pi^1_C\circ (\Delta_C +(-1)^{r_0}\Gamma_{\sigma_0})\circ  (\Delta_C +(-1)^{r_1}\Gamma_{\sigma_1})\circ \cdots\circ (\Delta_C +(-1)^{r_4}\Gamma_{\sigma_4})&\cr\
    \ \ \in A^1(C\times C)\ &.\cr\end{split}\end{equation}
    Each summand in \eqref{this} with more than one $r_i$ equal to 0 vanishes (because the quotient $C/\langle \sigma_i,\sigma_j\rangle$ is a rational curve which has no $h^1$ part in its motive). 
    Also, the summand in \eqref{this} with all $r_i$ equal to 1 vanishes (because $\sigma_0\circ\sigma_1\circ\cdots\sigma_4$ is the identity, which does not have an anti-invariant part). It follows that \eqref{this} simplifies to a 5-term decomposition
    \[ \pi^1_C= \pi^1_{C,(0)}+ \pi^1_{C,(1)}+\cdots+  \pi^1_{C,(4)}   \ \ \ \hbox{in}\ A^1(C\times C)\ ,  \]
    where $\pi^1_{C,(i)}$ is defined as the summand in \eqref{this} with $r_i=0$ (and hence all $r_j$, $j\not=i$ equal to 1).
 
 The elliptic curves $C_i$ have a decomposition $h(C_i)=h^0(C_i)\oplus h^1(C_i)\oplus h^2(C_i)$, induced by the degree 1 0-cycle $(\rho_i)_\ast(z)\in A^1(C_i)$. 
 For $k=0,2$, there are isomorphisms
      \begin{equation}\label{02} {}^t \Gamma_{\rho_i}\colon\  h^k(C_i)\cong h^k(C)\ \ \ \hbox{in}\ \MM_{\rm rat}\ \ \ (i=0,\ldots,4)\ .\end{equation}
      The above construction gives an isomorphism
      \begin{equation}\label{1}  {}^t \Gamma_{\rho_i}\colon\  h^1(C_i)\cong h^1_{(i)}(C)\ \ \ \hbox{in}\ \MM_{\rm rat}\ \ \ (i=0,\ldots,4)\ \end{equation}  
      (with inverse given by ${1\over 2} \Gamma_{\rho_i}$).
      Combining \eqref{02} and \eqref{1} proves the ``moreover'' part.
         \end{proof}
    
 \begin{remark}\label{+-} As in \cite{45}, let us write $h(C)^{+++++}$, $h(C)^{++++-}$ etc. for the submotive of $h(C)$ that is $\sigma_i$-invariant for $i=0,\ldots,4$, resp. that is $\sigma_4$-anti-invariant and $\sigma_i$-invariant for $i=0,\ldots,3$. The argument of Proposition \ref{splitting} then gives the decomposition
   \[ \begin{split} h(C)=&h(C)^{+++++}\oplus \\  &\ h(C)^{+----}\oplus h(C)^{-+---}\oplus h(C)^{--+--}\oplus  h(C)^{---+-}\oplus h(C)^{----+}\ ,\\
   \end{split}\]
   where $h^0(C)\oplus h^2(C)$ is equal to the first summand, and $h^1_{(0)}(C), \ldots, h^1_{(4)}(C)$ are equal to the other summands.
 \end{remark}

 \subsection{MCK decompositions}
\label{ss:mck}

\begin{definition}[Murre \cite{Mur}] Let $X$ be a smooth projective variety of dimension $n$. We say that $X$ has a {\em CK decomposition\/} if there exists a decomposition of the diagonal
   \[ \Delta_X= \pi^0_X+ \pi^1_X+\cdots +\pi_X^{2n}\ \ \ \hbox{in}\ A^n(X\times X)\ ,\]
  such that the $\pi^i_X$ are mutually orthogonal idempotents and $(\pi_X^i)_\ast H^\ast(X,\QQ)= H^i(X,\QQ)$.
  
  (NB: ``CK decomposition'' is shorthand for ``Chow--K\"unneth decomposition''.)
\end{definition}

\begin{remark} The existence of a CK decomposition for any smooth projective variety is part of Murre's conjectures \cite{Mur}, \cite{J4}. 
\end{remark}

\begin{definition}[Shen--Vial \cite{SV}] Let $X$ be a smooth projective variety of dimension $n$. Let $\Delta_X^{sm}\in A^{2n}(X\times X\times X)$ be the class of the small diagonal
  \[ \Delta_X^{sm}:=\bigl\{ (x,x,x)\ \vert\ x\in X\bigr\}\ \subset\ X\times X\times X\ .\]
  An {\em MCK decomposition\/} is a CK decomposition $\{\pi_X^i\}$ of $X$ that is {\em multiplicative\/}, i.e. it satisfies
  \begin{equation}\label{vani} \pi_X^k\circ \Delta_X^{sm}\circ (\pi_X^i\times \pi_X^j)=0\ \ \ \hbox{in}\ A^{2n}(X\times X\times X)\ \ \ \hbox{for\ all\ }i+j\not=k\ .\end{equation}
  Here $\pi^i_X\times \pi^j_X$ is by definition $(p_{13})^\ast(\pi^i_X)\cdot (p_{24})^\ast(\pi^j_X)\in A^{2n}(X^4)$, where $p_{rs}\colon X^4\to X^2$ denotes projection on $r$th and $s$th factors.
  
 (NB: ``MCK decomposition'' is shorthand for ``multiplicative Chow--K\"unneth decomposition''.) 
  
% A {\em weak MCK decomposition\/} is a CK decomposition $\{\pi^X_i\}$ of $X$ that satisfies
%    \[ \Bigl(\pi^X_k\circ \Delta_X^{sm}\circ (\pi^X_i\times \pi^X_j)\Bigr){}_\ast (a\times b)=0 \ \ \ \hbox{for\ all\ } a,b\in\ A^\ast(X)\ .\]
  \end{definition}
  
  \begin{remark}\label{rem:mck} Note that the vanishing \eqref{vani} is always true modulo homological equivalence; this is because the cup product in cohomology respects the grading.
  
    The small diagonal (seen as a correspondence from $X\times X$ to $X$) induces the {\em multiplication morphism\/}
    \[ \Delta_X^{sm}\colon\ \  h(X)\otimes h(X)\ \to\ h(X)\ \ \ \hbox{in}\ \MM_{\rm rat}\ .\]
 Suppose $X$ has a CK decomposition
  \[ h(X)=\bigoplus_{i=0}^{2n} h^i(X)\ \ \ \hbox{in}\ \MM_{\rm rat}\ .\]
  By definition, this decomposition is multiplicative if for any $i,j$ the composition
  \[ h^i(X)\otimes h^j(X)\ \to\ h(X)\otimes h(X)\ \xrightarrow{\Delta_X^{sm}}\ h(X)\ \ \ \hbox{in}\ \MM_{\rm rat}\]
  factors through $h^{i+j}(X)$.
  
  If $X$ has an MCK decomposition, then setting
    \[ A^i_{(j)}(X):= (\pi_X^{2i-j})_\ast A^i(X) \ ,\]
    one obtains a bigraded ring structure on the Chow ring: that is, the intersection product sends $A^i_{(j)}(X)\otimes A^{i^\prime}_{(j^\prime)}(X) $ to  $A^{i+i^\prime}_{(j+j^\prime)}(X)$.
    
      It is expected that for any $X$ with an MCK decomposition, one has
    \[ A^i_{(j)}(X)\stackrel{??}{=}0\ \ \ \hbox{for}\ j<0\ ,\ \ \ A^i_{(0)}(X)\cap A^i_{hom}(X)\stackrel{??}{=}0\ ;\]
    this is related to Murre's conjectures B and D, that have been formulated for any CK decomposition \cite{Mur}.

  The property of having an MCK decomposition is severely restrictive, and is closely related to Beauville's ``splitting property' conjecture'' \cite{Beau3}. 
  To give some idea: hyperelliptic curves have an MCK decomposition \cite[Example 8.16]{SV}, but the very general curve of genus $\ge 3$ does not have an MCK decomposition \cite[Example 2.3]{FLV2}. As for surfaces: a smooth quartic in $\PP^3$ has an MCK decomposition, but a very general surface of degree $ \ge 7$ in $\PP^3$ should not have an MCK decomposition \cite[Proposition 3.4]{FLV2}.
For a more detailed discussion, as well as examples of varieties with an MCK decomposition, we refer to \cite[Section 8]{SV}, as well as \cite{V6}, \cite{SV2}, \cite{FTV}, \cite{37}, \cite{38}, \cite{39}, \cite{40}, \cite{44}, \cite{45}, \cite{46}, \cite{48}, \cite{55}, \cite{FLV2}, \cite{g8}, \cite{59}, \cite{60}.
   \end{remark}

\subsection{MCK decomposition for curves}
\label{ssc}

\begin{proposition}\label{curveprop} Let $C$ be a curve of genus $g$, and $z\in A^1(C)$ a zero-cycle of degree 1. The following are equivalent:

\begin{enumerate}

\item The CK decomposition
  \[ \pi^0_C= z\times C\ ,\ \ \pi^2_C=C\times z\ ,\ \ \pi^1_C=\Delta_C-\pi^0_C-\pi^2_C \]
  is MCK;
  
  \item The modified small diagonal
       \[\begin{split} \Gamma_{3}(C,
			z):=&\Delta^{sm}_{C}-p_{12}^{*}(\Delta_{X})p_{3}^{*}(z)-p_{23}^{*}(\Delta_{C})p_{1}^{*}(z)-p_{13}^{*}(\Delta_{C})p_{2}^{*}(z)\\
			&+p_{1}^{*}(z)p_{2}^{*}(z)+p_{1}^{*}(z)p_{3}^{*}(z)+p_{2}^{*}(z)p_{3}^{*}(z) \ \ \ \ \ \hbox{in}\ A^2(C^3)\ \\
			\end{split}\]  
			is zero;
			
\item the class $[C]\in A^{g-1}(\jac(C))$ is in $A^{g-1}_{(0)}(\jac(C))$, where $C$ is embedded in the Jacobian $\jac(C)$ using the 0-cycle $z$, and $A^\ast_{(\ast)}(\jac(C))$ refers to the Beauville decomposition \cite{Beau}.
\end{enumerate}
\end{proposition}

\begin{proof} This is \cite[Proposition 3.1]{FLV2}.
\end{proof}

This leads to non-existence results:

\begin{corollary} Let $C$ be a very general curve of genus $>2$. Then $C$ does not admit an MCK decomposition.
\end{corollary}

\begin{proof} If $C$ admits an MCK decomposition, then Proposition \ref{curveprop} implies that the {\em Ceresa cycle\/}
    \[   [C] - (-1_{\jac(C)})_\ast [C]     \ \   \in\ A_1(\jac(C)) \]
must vanish. But Ceresa has shown that the Ceresa cycle of a very general curve of genus $>2$ is non-torsion modulo algebraic equivalence \cite{Cere}.
\end{proof}

  \begin{remark}\label{!} 
  In case $Y\subset\PP^{m}$ is a smooth {\em hypersurface\/} (of any degree), there is equality
    \begin{equation}\label{form} \Delta_Y\cdot (p_i)^\ast(h)= {\displaystyle\sum} a_k\, (p_1)^\ast (h^k)\cdot (p_2)^\ast (h^{m-k})\ \ \ \hbox{in}\ A^{m}(Y\times Y)\ ,\end{equation}
    with $a_k\in\QQ$,
as follows from the excess intersection formula. On the other hand, in case $Y\subset\PP^m$ is a complete intersection of codimension at least 2,
in general there is {\em no equality\/} of the form \eqref{form}. Indeed, let $C$ be a very general curve of genus $g\ge 4$. The Faber--Pandharipande cycle
  \[  FP(C):= \Delta_C\cdot (p_j)^\ast(K_C) - {1\over 2g-2} K_C\times K_C\ \ \ \in A^2(C\times C)\ \ \ \ \ (j=1,2) \]
  is homologically trivial but non-zero in $A^2(C\times C)$ \cite{GG}, \cite{Yin0} (this cycle $FP(C)$ is the ``interesting 0-cycle'' in the title of \cite{GG}). In particular, for the very general complete intersection $Y\subset\PP^3$ of bidegree $(2,3)$, the cycle
  \[ FP(Y):=  \Delta_Y\cdot (p_j)^\ast(h) - {1\over 6} h \times h\ \ \ \in A^2(Y\times Y) \]
  is homologically trivial but non-zero, and so there cannot exist an equality of the form \eqref{form} for $Y$.
  
  This is intimately related to MCK decompositions. Indeed, if the curve $Y$ had an MCK decomposition which is generically defined, then $K_Y\in A^1_{(0)}(Y)$ and hence
  $FP(Y)\in A^2_{(0)}(Y\times Y)\cong\QQ$, and so
  $FP(Y)$ would be zero.
  \end{remark}

\section{The main result}

\begin{theorem}\label{main} Any Humbert curve has an MCK decomposition.
\end{theorem}

\begin{proof} We consider the decomposition
 % \[ \pi^0_C\ ,\ \ \pi^1_C=\pi^1_{C,(0)}+\cdots+\pi^1_{C,(4)}\ ,\ \ \pi^2_C \]
 \[ h^0(C)\ ,\ \   h^2(C)\ ,\ \ h^1(C)=h^1_{(0)}(C)\oplus\cdots\oplus h^1_{(4)}(C)\  \]
  of Proposition \ref{splitting}, and we are going to verify this is an MCK decomposition.
  
  By definition, this means the following: given three motives
   \[  M_1, M_2, M_3\ \ \in \ \Bigl\{ h^0(C), h^2(C), h^1_{(0)}(C),\ldots, h^1_{(4)}(C)\Bigr\}\ \]
    of degree say $m_1$ resp. $m_2$ resp. $m_3\not=m_1+m_2$, we need to ascertain that the composition of the  ``multiplication map'' with the projection to $M_3$
       \[  \ M_1\otimes M_2\ \xrightarrow{\Delta^{sm}_C}\ h(C)\ \to\ M_3\]  
  is zero.
     
   As a first step, let us assume that there exists one $i$ such that
   \[  M_1, M_2, M_3\ \ \in \Bigl\{  h^0(C), h^2(C), h^1_{(i)}(C)\Bigr\} \ . \]
 In this case, we have a commutative diagram
        \[ \begin{array}[c]{ccccc}
                           M_1\otimes M_2 &   \xrightarrow{\Delta_{C}^{sm}} & h(C)  & \to& M_3 \\
                                   && && \\
                                 \ \  \ \ \ \ \  \ \ \ \uparrow{\scriptstyle {}^t \Gamma_{\rho_i}\times {}^t \Gamma_{\rho_i}}&&\ \ \ \uparrow{\scriptstyle {}^t \Gamma_{\rho_i}}\ && \ \ \ \downarrow {\scriptstyle \Gamma_{\rho_i}}\                                 \\
                                   && && \\
                                   h^{m_1}_{}(C_i)\otimes h^{m_2}_{}(C_i) & \xrightarrow{\Delta^{sm}_{C_i}} & h(C_i) & \to& h^{m_3}(C_i)  \\
                 \end{array}\]  
  where $C_i=C/\langle \sigma_i\rangle$ is an elliptic curve. The left square commutes because of contravariance of the intersection product; the right square commutes because of Lieberman's lemma and the fact that $(\rho_i\times\rho_i)_\ast (p)=\pi^{m_3}_{C_i}$, where $p$ is the projector on $M_3\in \{  h^0(C), h^2(C), h^1_{(i)}(C)\} $. Now the decomposition $h(C_i)=h^0(C_i)\oplus h^1(C_i)\oplus h^2(C_i)$ of Proposition \ref{splitting} is MCK (indeed, the 0-cycle used to define this decomposition is supported on the Weierstrass points), and so the composition of the two lower horizontal arrows is zero (since by assumption $m_3\not=m_1+m_2$). But the outer two vertical arrows are isomorphisms (Proposition \ref{splitting}), and so
the composition of the top horizontal arrows is zero as well.   

 It remains to consider all the cases where there are at least two different $h^1_{(i)}(C)$ among $M_1, M_2, M_3$. For notational ease, we will assume
   \[  h^1_{(0)}(C), h^1_{(1)}(C)\ \ \in\ \{ M_1, M_2,M_3\}\ .\]
 
 Let us assume $M_1= h^1_{(0)}(C)$ and $M_2= h^1_{(1)}(C)$.  
  In terms of the notation of Remark \ref{+-}, the ``multiplication map''
       \[ \Delta_{C}^{sm}\colon\ \  h^1_{(0)}(C)\otimes h^1_{(1)}(C)= h(C)^{+----}\otimes h(C)^{-+---}  \  \to\ h(C)\]
     factors over $h^{--+++}(C)$, which is zero (cf. Remark \ref{+-}).
     
 Next, let us assume $M_1= M_2=h^1_{(0)}(C)$ (and so $M_3= h^1_{(1)}(C)$). Then
    \[ \Delta_{C}^{sm}\colon\ \  h^1_{(0)}(C)\otimes h^1_{(0)}(C)= h(C)^{+----}\otimes h(C)^{+----}  \  \to\ h(C)\]
     factors over $h^{+++++}(C)$, and so
     \[ M_1\otimes M_2\ \xrightarrow{\Delta_C^{sm}}\ h(C)\ \to\  h(C)^{-+---}=M_3 \]
     is zero.
     
  Finally, let us assume $M_1= h^1_{(0)}(C)$, $M_3= h^1_{(1)}(C)$, and $M_2\in\{  h^0(C), h^2(C)\}$.
  Then $M_1=h(C)^{+----}$ and
  $M_2$ is a submotive of $h(C)^{+++++}$, and so
  \[        M_1\otimes M_2\ \xrightarrow{\Delta_C^{sm}}\ h(C) \]
  factors over $h(C)^{+----}$, and we conclude that composing with the projection to $M_3=h(C)^{-+---}$ yields zero.
      This ends the proof.
     \end{proof}

 \begin{theorem}\label{main1} Let $C$ be a Humbert curve, and $m\in\NN$. Let
  \[ R^\ast(C^m):= \bigl\langle  (p_{ij})^\ast(\Delta_C), (p_k)^\ast(K_C)\bigr\rangle \ \ \subset\ A^\ast(C^m) \]
  be the $\QQ$-subalgebra generated by (pullbacks of) the diagonal $\Delta_C\subset C\times C$ and (pullbacks of) the canonical divisor $K_C$. The cycle class map induces injections
   \[ R^\ast(C^m)\ \hookrightarrow\ H^\ast(C^m,\QQ) \]
   for all $m\in\NN$.
  \end{theorem}
  
  \begin{proof} The result is actually true for any curve with an MCK decomposition; this is the content of \cite[Proposition 3.5]{FLV2}. Let us briefly recap the argument, which is essentially due to Tavakol \cite{Ta2}, \cite{Ta}, and formally similar to analogous results for cubic hypersurfaces \cite[Section 2.3]{FLV3} and for K3 surfaces \cite{Yin}.  
  
%  argument is shamelessly borrowed from the analogous result for cubic hypersurfaces \cite[Section 2.3]{FLV3}, which in turn is inspired by analogous results for hyperelliptic curves \cite{Ta2}, \cite{Ta} and for K3 surfaces \cite{Yin}.

As in \cite[Section 2.3]{FLV3}, let us write $o\in A^{1}(C)$ for the degree 1 0-cycle of Proposition \ref{splitting}, and
  \[ \tau:= \Delta_C -  o\times C - C\times o \ \ \in\ A^{1}(C\times C) \]
  (this cycle $\tau$ is just the projector $\pi^1_C$ on the motive $h^{1}(C)$ considered above).
Moreover, let us write 
  \[ \begin{split}   o_i&:= (p_i)^\ast(o)\ \ \in\ A^{1}(C^m)\ ,\\
                      %  h_i&:=(p_i)^\ast(h)\ \ \in \ A^1(Y^m)\ ,\\
                         \tau_{i,j}&:=(p_{ij})^\ast(\tau)\ \ \in\ A^{1}(C^m)\ .\\
                         \end{split}\]
We now define the $\QQ$-subalgebra
  \[ \bar{R}^\ast(C^m):=\bigl\langle o_i,  \tau_{i,j}\bigr\rangle\ \ \ \subset\ H^\ast(C^m,\QQ)\ \ \ \ \ (1\le i\le m,\ 1\le i<j\le m)\ ; \]
 this is the image of $R^\ast(C^m)$ in cohomology. One can prove (just as \cite[Lemma 2.12]{FLV3} and \cite[Lemma 2.3]{Yin}) that the $\QQ$-algebra $ \bar{R}^\ast(C^m)$
  is isomorphic to the free graded $\QQ$-algebra generated by $o_i,\tau_{i,j}$, modulo the following relations:
    \begin{equation}\label{E:X'}
			o_i\cdot o_i = 0,   \quad 
			%h_i^{2g-1} =4\,o_i\,;
			\end{equation}
			\begin{equation}\label{E:X2'}
			\tau_{i,j} \cdot o_i = o_i\cdot o_j , \quad \tau_{i,j} \cdot \tau_{i,j} = -b\, o_i\cdot o_j
			%, \quad \tau_{i,j}=\tau_{j,i}
			%e\, o_i\cdot o_j
			\,;
			\end{equation}
			\begin{equation}\label{E:X3'}
			\tau_{i,j} \cdot \tau_{i,k} = \tau_{j,k} \cdot o_i\,;
			\end{equation}
			\begin{equation}\label{E:X4'}
			\sum_{\sigma \in \mathfrak{S}_{b+2}} 
			%\hbox{sign}(\sigma) 
			\prod_{i=1}^{b/2+1} \tau_{\sigma(2i-1), \sigma(2i)} = 0\, ,
			\end{equation}
			where $b:= \dim H^{1}(C,\QQ)= 10$.
			% and $e:=2g-b$ is the topological euler characteristic of $Y$ which is 0 !!!!

To prove Theorem \ref{main1}, it suffices to check that these relations are verified modulo rational equivalence.
Relation \eqref{E:X'} is clearly true, for reasons of dimension.

The relations \eqref{E:X2'} are pullbacks of relations taking place in $R^\ast(C^2)$. 
The first relation of \eqref{E:X2'} concerns the Faber--Pandharipande cycle of Remark \ref{!}. To prove this relation modulo rational equivalence, we note that $\tau_{}$ is in 
$A^1_{(0)}(C^2)$ (where $C^2$ is given the product MCK decomposition), and $o$ is in $A^1_{(0)}(C)$. It follows that 
  \[  \tau_{} \cdot (p_1)^\ast(0)- (p_1)^\ast(o)\cdot (p_2)^\ast(o) \ \ \in A^2_{(0)}(C ^2)\ ,\]
  and it is readily checked that $A^2_{(0)}(C ^2)$ injects into cohomology. This argument also applies to the second relation of \eqref{E:X2'}. 
  %The last relation of  \eqref{E:X2'} is clear.
    %The first and the last relations are trivially verified, because ($Y$ being Fano) one has
%$A^{4g-2}(Y^2)=\QQ$. As for the second relation of \eqref{E:X2'}, this is relation \eqref{hyp}.
% follows from the Franchetta property for $Y\times Y$ (Proposition \ref{Fr}), or one can prove it directly
%(Alternatively, one can deduce the second relation from the MCK decomposition: the product $\tau_{} \cdot h_i$ lies in $A^4_{(0)}(Y^2)$, and it is readily checked that $A^4_{(0)}(Y^2)$ injects into $H^8(Y^2,\QQ)$.)
%we use that $Y$ (and hence $Y^2$) has an MCK decomposition.
%Since $\Delta_Y$ and $h$ are in the subring $A^\ast_{(0)}(Y^2)$, the product $\tau_{} \cdot h_i$ lies in $A^4_{(0)}(Y^2)$.
%But $A^4_{(0)}(Y^2)$
 %  injects into $H^8(Y^2,\QQ)$ (Lemma \ref{40} below), and so $\tau_{} \cdot h_i=0$ in $A^4(Y^2)$.  
   
   Relation \eqref{E:X3'} is the pullback of a relation taking place in $R^\ast(C^3)$ and also follows from the MCK decomposition. Indeed, we have
   \[  \Delta_C^{sm}\circ (\pi^{1}_C\times\pi^{1}_C)=   \pi^{2}_C\circ \Delta_C^{sm}\circ (\pi^{1}_C\times\pi^{1}_C)  \ \ \ \hbox{in}\ A^{2}(C^3)\ ,\]
   and (using Lieberman's lemma) this translates into
   \[ (\pi^{1}_C\times \pi^{1}_C\times\Delta_C)_\ast    \Delta_C^{sm}  =   ( \pi^{1}_C\times \pi^{1}_C\times\pi^{2}_C)_\ast \Delta_C^{sm}                
                                            \ \ \ \hbox{in}\ A^{2}(C^3)\ ,\]
   which means that
   \[  \tau_{1,3}\cdot \tau_{2,3}= \tau_{1,2}\cdot o_3\ \ \ \hbox{in}\ A^{2}(C^3)\ .\]
   
  Finally, relation \eqref{E:X4'}, which takes place in $R^\ast(C^{b+2})$
   follows from the Kimura finite dimensionality relation \cite{Kim}: relation \eqref{E:X4'} expresses the vanishing
    \[ \sym^{b+2} H^{1}(C,\QQ)=0\ ,\]
    where $H^{1}(C,\QQ)$ is considered as a super vector space.
 This relation is also verified modulo rational equivalence
 (i.e., relation \eqref{E:X4'} is also true in $A^{(b+2)}(C^{b+2})$): indeed, relation \eqref{E:X4'} involves a cycle in
   \[ A^\ast(\sym^{b+2} h^1(C))\ ,\]
  but ($\sym^{b+1} h^1(C)$ and so a fortiori) $\sym^{b+2} h^1(C)$ is $0$ 
 because curves are Kimura finite-dimensional \cite{Kim}.
  %More precisely, Kimura finite-dimensionality of $Y$ implies there is a relation
 %  \[  \sym^{b+1} \pi^{2g-1}_Y = \sum_{\sigma \in \mathfrak{S}_{b+1}} \prod_{i=1}^{b+1} \tau_{i, b+1+\sigma(i)} = 0 \ \ \ \hbox{in}\ A^{(b+1)(2g-1)}(Y^{b+1}\times Y^{b+1})\ ,\]
 %  and (because cup product is graded commutative) this is sent to the relation  \eqref{E:X4'} under the cycle class map. 
    %\cite[Theorem 4]{43}. 
 This ends the proof.
%  
% \medskip
% Above we have used the following lemma:
%  
% \begin{lemma}\label{40} The cycle class map induces an injection
%   \[ A^4_{(0)}(Y^2)\ \hookrightarrow\ H^8(Y^2,\QQ)\ .\] 
% \end{lemma} 
% 
% To prove the lemma, we note that by definition
%  \begin{equation}\label{proj} A^4_{(0)}(Y^2)= (\pi^6_Y\times \pi^2_Y + \pi^4_Y\times \pi^4_Y + \pi^2_Y\times\pi^6_Y)_\ast A^4(Y^2) \ .\end{equation}
%  These three direct summands land in different pieces of the K\"unneth decomposition of $H^8(Y^2,\QQ)$, and so we may check injectivity summand by summand. By construction, $\pi^{2j}_Y$ is supported on $V_j\times W_j$ where $V_j\subset Y$ is smooth of dimension $j$ and $W_j\subset Y$ is smooth of codimension $j$. Thus, the action of each projector in \eqref{proj} factors over $A^0(W_j\times W_{4-j})\cong\QQ$, which proves the lemma.
 \end{proof}

\section{Closing remarks}

\begin{remark} The argument of Theorem \ref{main}, using a motivic decomposition coming from a large automorphism group, is formally similar to that of \cite{45}, where an MCK decomposition was established for special Horikawa surfaces. In its turn, \cite{45} was inspired by the argument of \cite[Section 7.1]{FLV2}, where MCK was established for certain Todorov surfaces. In practice, this type of argument tends to be applicable over certain special loci in the moduli space (where the automorphism group gets large), and not for the general element of some locally complete family.
\end{remark}

\begin{remark} It seems likely that Humbert--Edge curves (of any type $n$) admit an MCK decomposition. However, the naive ``$--=+$'' argument proving Theorem \ref{main} does not suffice to settle this for $n>4$. For example, let $C$ be a Humbert--Edge curve of type 5 (the genus of $C$ is 17), with commuting involutions $\sigma_0,\ldots,\sigma_5$. Using \cite[Theorem 2.1]{AAR}, one can show there is a decomposition
  \[ h^1(C)=h(C)^{------} \oplus h(C)^{++----}\oplus \hbox{perm.}\ \ \ \hbox{in}\ \MM_{\rm rat}\ ,\]
  where ``perm.'' means all possible permutations of 2 plus signs and 4 minus signs. Now it is not clear whether the map
  \[ \Delta_C^{sm}\colon\ \ h(C)^{++----}\otimes h(C)^{--++--}\ \to\ h(C)^{----++} \]
  is the zero map (which is a necessary condition for the decomposition to be MCK).
  
  Some more sophisticated reasoning is called for here; perhaps it is possible to show that $C$ satisfies condition $(\ast)$ of \cite{FV} ?
\end{remark}

\begin{remark} In \cite{BLLS} and \cite{Bcurve2}, examples are given of non-hyperelliptic curves for which the Ceresa cycle is torsion {\em modulo algebraic equivalence\/}.
In view of \cite[Remark 3.4]{FLV2}, these curves hence have an MCK decomposition modulo algebraic equivalence (and the argument of Theorem \ref{main1} gives that the tautological ring modulo algebraic equivalence of these curves injects into cohomology).
It would be interesting to settle whether these results can be lifted to rational equivalence.
\end{remark}

\begin{remark}\label{new}(Added in revision) The recent preprint \cite{QZ} contains many examples of non-hyperelliptic curves with vanishing Ceresa cycle and vanishing modified small diagonal {\em modulo rational equivalence\/}. In view of Proposition \ref{curveprop}, these examples also have an MCK decomposition.
\end{remark}

\vskip0.5cm
\begin{nonumberingda} The author states that this is not applicable, as there are no associated data.
\end{nonumberingda}

\vskip0.5cm
\begin{nonumberingcon} The author states that there is no conflict of interest.
\end{nonumberingcon}

 \vskip0.5cm
\begin{nonumberingt} The author thanks Yoyo for a wonderful Japanese dinner on June 18 2021.
\end{nonumberingt}

\end{document}